\documentclass{amsart}
\usepackage{amsfonts}
\usepackage{amssymb}
\usepackage{mathrsfs}
\usepackage{amsmath}
\usepackage{cite}

\newtheorem{theorem}{Theorem}

\newtheorem{remark}[theorem]{Remark}

\newtheorem{corollary}[theorem]{Corollary}

\newtheorem{definition}[theorem]{Definition}

\newtheorem{lemma}[theorem]{Lemma}
\newtheorem{proposition}[theorem]{Proposition}

\numberwithin{equation}{section} \numberwithin{theorem}{section}

\renewcommand{\oddsidemargin}{5mm}

\begin{document}

\title
{The blow up analysis  \\of the general curve shortening flow}
\author{RongLi Huang$^{1,2}$}
\address{1. School of Mathematical Sciences, Beijing Normal University,
Laboratory of Mathematics and Complex Systems, Ministry of Education,
Beijing 100875, People's Republic of China}
\email{hrl602@mail.bnu.edu.cn}
\address{2. Institute of Mathematics, Fudan University,
Shanghai 200433, People's Republic of China}
\email{huangronglijane@yahoo.cn}
\author{JiGuang Bao$^3$}
\address{3. Corresponding author. School of Mathematical Sciences, Beijing Normal University,
Laboratory of Mathematics and Complex Systems, Ministry of Education,
Beijing 100875, People's Republic of China}
\email{jgbao@bnu.edu.cn }

\date{}
\begin{abstract}  It is shown that the curvature function satisfies a nonlinear evolution equation
under the general curve shortening flow and
 a detailed asymptotic behavior of the closed curves is presented
when they contract to a point in finite time.
\end{abstract}

\maketitle

\section{Introduction}
The curve shortening flows have been studied by many authors and have
many applications (cf. \cite{BS}). M.E.Gage and R.S.Hamilton  discussed
the mean curvature flow in one dimensional case (cf. \cite{GH}). The flow
is given by the equation
\begin{equation*}\label{1.01}
 v=k,
\end{equation*}
where $v$ and $k$ are, respectively, the normal velocity and inward
curvature of the plane curve. They proved that a convex closed curve
stays convex and smooth and shrinks to a point in finite time with
the limiting shape of a circle. In this paper, we generalized some results on
curve shortening flow (cf.\cite{GH}) to the  generalized curve shortening flow (cf. \cite{B})
\begin{equation}\label{e1.1}
 v=|k|^{p-1}k,
\end{equation}
where  $p$ is the positive  number. When $p=\frac{1}{3}$, the flow
(\ref{e1.1}) is the affine plane curve evolution (cf. \cite{COT}, \cite{ST}). More
generally, we will also study the non-homogeneous flow (cf. \cite{B}) which is
given by
\begin{equation}\label{e1.2}
 v=G(k)k,
\end{equation}
where $G(\cdot)$ is a positive function on $(0,\infty)$. Chia-Hsing
Nien and Dong-Ho Tsai had proved that the self-similar solutions
under the contraction flow could happen only as (\ref{e1.1}) (cf. \cite{NT}).
B.Andrews had studied some nonlinear expansion of contraction flow
and obtained the limiting self-similar solutions (cf. \cite{B}). But the
asymptotic behavior of solutions of the flow (\ref{e1.1}), (\ref{e1.2}) are
little known. The main objective of this paper is to analyze the
asymptotic behavior of the curvature under the generalized curve shortening  flow .

Let $\mathbb{S}^{1}$ be an unit circle in the plane, and
\begin{equation*}
 \gamma_{0}: \mathbb{S}^{1}\rightarrow \mathbb{R}^{2} ,
\end{equation*}
be closed convex curve in the plane. We look for a  family of closed curves
\begin{equation*}
 \gamma(u,t): \mathbb{S}^{1}\times [0,T)\rightarrow \mathbb{R}^{2} ,
\end{equation*}
which satisfies
\begin{equation*}
\left\{ \begin{aligned} \frac{\partial \gamma}{\partial t}(u,
t)&=|k|^{p-1}kN , \quad
&u\in \mathbb{S}^{1}, \quad &t\in[0, T),\\
\gamma(u,0)&=\gamma_{0}(u), \quad &u\in \mathbb{S}^{1},\quad &t=0,
\end{aligned} \right.
\end{equation*}
where $p$ is a positive number, $k(\cdot ,t)$ is the inward curvature of the
plane curve $\gamma(\cdot ,t)$ and $N(\cdot ,t)$ is the  unit inward normal
vector. More generally, we consider $\gamma(\cdot ,t)$ satisfying
\begin{equation}\label{e1.3}
\left\{ \begin{aligned} \frac{\partial \gamma}{\partial t}(u,
t)&=G(k)kN , \quad
&u\in \mathbb{S}^{1},\quad &t\in[0, T),\\
\gamma(u,0)&=\gamma_{0}(u), \quad &u\in \mathbb{S}^{1}, \quad &t=0,
\end{aligned} \right.
\end{equation}
where $G$ is a  positive, non-decreasing smooth function on
$(0,\infty)$.

In the following sections we assume that $A(t)$ is  the area
of a bounded domain enclosed by the curve $\gamma(\cdot ,t)$,
$L(t)$ is  the length of $\gamma(\cdot, t)$,  $r_{out}(t)$ and $r_{in}(t)$
are respectively the radii of the largest circumscribed circle and the smallest circumscribed circle
of $\gamma(\cdot ,t)$. Define
$$ k_{\max}(t)=\max\{k(u,t)\mid u\in \mathbb{S}^{1}\}, $$
$$ k_{\min}(t)=\min\{k(u,t)\mid u\in \mathbb{S}^{1}\}. $$

Firstly we introduce the existence theorem, which belongs to Ben.Andrews
(cf. Theorem $\Pi$4.1, Proposition $\Pi$ 4.4 in \cite{B}).
\begin{proposition}\label{p1.1}
Let $\gamma_{0}$ be a closed strictly convex curve. Then the unique
classical  solution $\gamma(\cdot,t)$ of (\ref{e1.3}) exists only at finite time interval $[0,\omega)$, and
 the solution $\gamma(\cdot,t)$ converges to a point $\vartheta$ as
$t\rightarrow \omega$ and $A(t)$, $k_{\max}(t)$  satisfy the
following properties:
$$\forall t\in [0,\omega), A(t)>0, k_{\max}(t)<+\infty, $$
$$\lim_{t\rightarrow \omega}A(t)=0, \ \ \ \ \lim_{t\rightarrow \omega}k_{\max}(t)=+\infty.$$
\item As $t\rightarrow \omega$, the normalized curves
$$ \eta(\cdot,t)=\sqrt{\frac{\pi}{A(t)}}\gamma(\cdot,t) $$
converges to the unit circle centered at the point $\vartheta$.

\end{proposition}
In this paper we always assume that the initial curve satisfies
the conditions of Proposition \ref{p1.1}, and that $G(x)$ is a function on $(0,\infty)$ satisfying

{\bf ($H1$)}\ $G(x) \in C^{3}(0,\infty)$, $G'(x)\geq 0$ and $G(x)>0$ for $x\in (0,\infty)$.

{\bf ($H2$)}\ $G(x)x^{2}$ is convex in $(0,\infty)$ and there is a positive constant $C_{0}$ such that
\begin{equation*}
G'(x)x\leq C_{0}G(x),\,\,\,\mathrm{for}\,\,   \mathrm{sufficiently} \,\,\mathrm{large} \,\,x.
\end{equation*}
We now state the main theorem of this paper.
\begin{theorem}\label{t1.1}
Suppose G(x) satisfies $(H1)$ and $(H2)$.  Let
$\gamma(\cdot,t)$ be the solution for Proposition \ref{p1.1}. Then the
following hold:
\begin{enumerate}
\item
$\displaystyle\lim_{t\rightarrow\omega}\frac{r_{in}(t)}{r_{out}(t)}=1.$

\item
$\displaystyle\lim_{t\rightarrow\omega}\frac{k_{\min}(t)}{k_{\max}(t)}=1.$

\item $\displaystyle\lim_{t\rightarrow \omega}\frac{1}{\omega-t}\int_{k(\theta,t)}^{+\infty}\frac{dx}{G(x)x^{3}}=1 $
is uniformly convergent on $\mathbb{S}^{1} $.
\end{enumerate}
\end{theorem}
\begin{remark}
Let $G(x)=| x|^{p-1}$ with $p\geq 1$ in Proposition 1.1.  Then
\begin{equation}\label{e1.4}
k(\theta,t)[(p+1)(\omega-t)]^{\frac{1}{p+1}}\,\, \mathrm{converges} \,\, \mathrm{uniformly }\,\,
\mathrm{to}\,\, 1 \,\, \mathrm{as}\,\, t\rightarrow \omega.
\end{equation}
for uniformly $\theta$ in $\mathbb{S}^{1} $.
\end{remark}
\begin{remark}\label{de 2.02}
When $p=1,$  by (\ref{e1.4}) it shows that the asymptotic formula about curvature
function of the curve shortening flow (\ref{e1.1}) which was discovered
firstly by M.Gage and R.S.Hamilton (cf. Corollary 5.6 in \cite{GH}).
\end{remark}
This paper is organized as follows: In the next section we
transfer the flow (\ref{e1.3}) into an initial PDEs problem and establish
some monotone geometric inequality. Section 3 is devoted to the main
contribution of our article, that is, the global Harnack inequality
of the curvature function according to the flow (\ref{e1.3}), see Lemma
3.6. And then we complete the proof of Theorem \ref{t1.1} by making use of
Gage-Hamilton's methods (cf. \cite{GH}).

\section{Evolutions}

Using  the idea in \cite{GH}, we can drive the evolution equations under the flow (\ref{e1.3}) for the length
and the curvature of the curves, and the area enclosed by the curves.

Let the curve be $\gamma(u)=(x(u),y(u))$ with parameter $u$ (modulo
2$\pi$) and $s$ be an arc-length parameter along the curve
$\gamma(u)$ which is unique up to a constant. Then
\begin{equation*}
ds=vdu,\quad \frac{\partial}{\partial s}=\frac{1}{v}\frac{\partial}{\partial u},
\end{equation*}
where $$v=\sqrt{(\frac{\partial x}{\partial u})^{2}+(\frac{\partial y}{\partial u})^{2}}.$$

Suppose $\xi$ and $N$ are the unit tangent vector and the unit inward normal vector of the curve.
Then the Frenet equations (cf. \cite{Ca}) are
\begin{equation*}
\frac{\partial \xi}{\partial u}=vkN , \ \ \frac{\partial N}{\partial u}=-vk\xi,
\end{equation*}
or
\begin{equation*}
\frac{\partial \xi}{\partial s}=kN , \ \ \frac{\partial N}{\partial s}=-k\xi,
\end{equation*}
where $k(u)$ is  the inward curvature of $\gamma(u)$.

Let $\theta$ be the tangent angle of the curve $\gamma(u)$ to the $x$-axis. We drive the following useful
formula.

\begin{lemma}\label{l2.1}
The geometric quantities $v$, $L$, $\zeta$, $N$, $\theta$ and $k$ of
the flow (\ref{e1.3}) evolve according to
\begin{enumerate}
\item
$\displaystyle\frac{\partial v}{\partial t}=-G(k)k^{2}v$,

\item
$\displaystyle\frac{d L}{d t}=-\int_{0}^{L}G(k)k^{2}ds $,

\item
$\displaystyle\frac{\partial }{\partial t}\frac{\partial }{\partial s}=\frac{\partial }{\partial t}\frac{\partial }{\partial s}+
G(k)k^{2}\frac{\partial }{\partial s} $,

\item
$\displaystyle\frac{\partial\xi}{\partial t}=(G(k)k)'\frac{\partial k}{\partial
s}N,\ \frac{\partial N}{\partial t}=-(G(k)k)'\frac{\partial
k}{\partial s}\xi $,

\item
$\displaystyle\frac{\partial\theta}{\partial t}=(G(k)k)'\frac{\partial k}{\partial
s},\ \frac{\partial\theta}{\partial s}=k$, where
$\xi=(\cos\theta,\sin\theta)$,

\item
$\displaystyle\frac{dA}{dt}=-\int^{L}_{0}G(k)kds,$

\item
$\displaystyle\frac{\partial k}{\partial t}=\frac{\partial}{\partial
s}\left((G(k)k)'\frac{\partial k}{\partial s}\right)+G(k)k^{3}.$
\end{enumerate}
\end{lemma}
\begin{proof}
Let $\langle\,\cdot,\cdot\rangle$ be the inner product in $\mathbb{R}^{2}$.

i) \ By (\ref{e1.3}) and the Frenet equations we have
\begin{equation*}\label{1.03}\aligned
\frac{\partial}{\partial t}(v^{2})
&=\frac{\partial}{\partial t}\langle\frac{\partial\gamma}{\partial\omega},\frac{\partial\gamma}{\partial\omega}\rangle \\
&=2\langle\frac{\partial\gamma}{\partial\omega},\frac{\partial^{2}\gamma}{\partial \omega\partial t}\rangle \\
&=2\langle v\xi, \frac{\partial}{\partial \omega}(G(k)kN)\rangle\\
&=2\langle v\xi, \frac{\partial}{\partial\omega}(G(k)k)N+G(k)k\frac{\partial N}{\partial\omega}\rangle \\
&=2\langle v\xi, \frac{\partial}{\partial\omega}(G(k)k)N-G(k)k^{2}v\xi\rangle\\
&=-2G(k)k^{2}v^{2}. \endaligned
\end{equation*}
This implies that  the  identity (i) holds.

ii) \ Since $L=\int^{2\pi}_{0}vd\omega$, then from (i) there holds
\begin{equation*}\label{1.03}
\frac{\partial L}{\partial t}=\int^{2\pi}_{0}\frac{\partial
v}{\partial
t}d\omega=-\int^{2\pi}_{0}G(k)k^{2}vd\omega=-\int^{L}_{0}G(k)k^{2}ds.
\end{equation*}

iii) \ By (i) we get
\begin{equation*}\label{1.03} \aligned
\frac{\partial}{\partial t}\frac{\partial}{\partial s}
&=\frac{\partial}{\partial t}(\frac{1}{v}\frac{\partial}{\partial \omega}) \\
&=-\frac{1}{v^{2}}\frac{\partial v}{\partial t}\frac{\partial}{\partial \omega}
+\frac{1}{v}\frac{\partial }{\partial t}\frac{\partial}{\partial \omega}\\
&=-\frac{1}{v^{2}}(-G(k)k^{2}v)\frac{\partial}{\partial\omega}+\frac{\partial }{\partial s}\frac{\partial}{\partial t}\\
&=G(k)k^{2}\frac{\partial}{\partial s}+\frac{\partial}{\partial s}\frac{\partial}{\partial t}
.\qquad\qquad\qquad\qquad\qquad\qquad\qquad\qquad\quad \endaligned
\end{equation*}

iv) \ It follows from (\ref{e1.3}) and (iii) that
$$\frac{\partial \xi}{\partial
t}=\frac{\partial}{\partial t}\frac{\partial}{\partial
s}\gamma=\frac{\partial}{\partial s}\frac{\partial}{\partial
t}\gamma+G(k)k^{2}\frac{\partial}{\partial
s}\gamma=\frac{\partial}{\partial s}(G(k)kN)+G(k)k^{2}\xi. $$ Then
by Frenet equations we have
\begin{equation}\label{1.03} \aligned
\frac{\partial \xi}{\partial t}&=\frac{\partial}{\partial
s}(G(k)k)N+G(k)k\frac{\partial N}{\partial s}+G(k)k^{2}\xi
\qquad\qquad\qquad\quad\quad \\
&= (G(k)k)'\frac{\partial k}{\partial
s}N-G(k)k^{2}\xi+G(k)k^{2}\xi\\
&=(G(k)k)'\frac{\partial k}{\partial s}N . \endaligned
\end{equation}
In terms of $\langle\xi,N\rangle\equiv 0$ we obtain
\begin{equation*}\label{1.03}
0\equiv\frac{\partial}{\partial t}\langle\xi,N\rangle=\langle\frac{\partial
\xi}{\partial t},N\rangle+\langle\xi,\frac{\partial N}{\partial
t}\rangle.
\end{equation*}
Thus
$$0=\langle(G(k)k)'\frac{\partial k}{\partial s}N,N\rangle+\langle\xi,\frac{\partial
N}{\partial t}\rangle=(G(k)k)'\frac{\partial k}{\partial
s}+\langle\xi,\frac{\partial N}{\partial t}\rangle .$$
From $0\equiv\langle\displaystyle\frac{\partial N}{\partial t}, N\rangle$ it
follows that there exists $\lambda$ such that
$$\frac{\partial N}{\partial t}=\lambda\xi,$$
and then combining with the above equality we have
$$\lambda=-(G(k)k)'\frac{\partial k}{\partial s},$$
and
\begin{equation*}\label{1.03} \frac{\partial N}{\partial
t}=-(G(k)k)'\frac{\partial k}{\partial s}\xi.
\end{equation*}

v) \ Since $\xi=(\cos\theta, \sin\theta)$, then
$N=(-\sin\theta,\cos\theta)$,  we obtain
$$\frac{\partial\xi}{\partial
t}=(-\sin\theta,\cos\theta)\frac{\partial\theta}{\partial
t}=\frac{\partial\theta}{\partial t}N $$ and comparing it with (iv)
we conclude that
\begin{equation*}\label{1.03} \frac{\partial\theta}{\partial
t}=(G(k)k)'\frac{\partial k}{\partial s}.
\end{equation*}
In other cases, $$\frac{\partial\xi}{\partial
s}=(-\sin\theta,\cos\theta)\frac{\partial\theta}{\partial
s}=\frac{\partial\theta}{\partial s}N. $$ It follows from
$\displaystyle\frac{\partial\xi}{\partial s}=kN$ that
\begin{equation*}\label{1.03}
\frac{\partial\theta}{\partial s}=k .
\end{equation*}

vi) \ Consider the closed curve
$\gamma=\{(x(u),y(u))|u\in \mathbb{S}^{1}\}$ in
$\mathbb{R}^{2}$. Then it is well known that the area of the domain
by the curve $\gamma$ can be expressed by the formula
\begin{equation*}
A=\frac{1}{2}\int_{0}^{2\pi}(x\frac{\partial y}{\partial
u}-y\frac{\partial x}{\partial
u})du=-\frac{1}{2}\int_{0}^{2\pi}\langle\gamma
(u),vN\rangle du .
\end{equation*}
Then
\begin{equation*}
\frac{dA}{dt}=-\frac{1}{2}\int_{0}^{2\pi}\langle\frac{\partial\gamma}{\partial
t} ,vN\rangle du -\frac{1}{2}\int_{0}^{2\pi}\langle\gamma
,\frac{\partial v}{\partial t}N\rangle du
-\frac{1}{2}\int_{0}^{2\pi}\langle\gamma ,\frac{\partial N}{\partial
t}v\rangle du.
\end{equation*}
By (\ref{e1.3}) and (i), (iv) we obtain
\begin{equation*}\aligned
\frac{dA}{dt}&=-\frac{1}{2}\int_{0}^{2\pi}\langle G(k)kN ,vN\rangle du
+\frac{1}{2}\int_{0}^{2\pi}\langle\gamma ,G(k)k^{2}vN\rangle du
+\frac{1}{2}\int_{0}^{2\pi}\langle\gamma,\frac{\partial}{\partial s}(G(k)k)v\xi\rangle du\\
&=-\frac{1}{2}\int_{0}^{2\pi}G(k)kvdu
+\frac{1}{2}\int_{0}^{2\pi}\langle\gamma ,G(k)k^{2}N\rangle vdu
+\frac{1}{2}\int_{0}^{2\pi}\langle\gamma,\frac{\partial}{\partial u}(G(k)k)\xi\rangle du\quad\\
&=-\frac{1}{2}\int_{0}^{2\pi}G(k)kvdu
+\frac{1}{2}\int_{0}^{2\pi}\langle\gamma ,G(k)k^{2}N\rangle vdu \\
& \ \ \ \ \ \ \ \ \ +\frac{1}{2}\int_{0}^{2\pi}\langle\gamma,\frac{\partial}{\partial u}(G(k)k\xi)\rangle du
-\frac{1}{2}\int_{0}^{2\pi}\langle\gamma ,G(k)k\frac{\partial\xi}{\partial u}\rangle du.
\endaligned\end{equation*}
Thus from $\displaystyle\frac{\partial \xi}{\partial w}=vkN$ we
arrive at
\begin{equation*}\aligned
\frac{dA}{dt}&=-\frac{1}{2}\int_{0}^{2\pi}G(k)kvdu
+\frac{1}{2}\int_{0}^{2\pi}\langle\gamma ,G(k)k^{2}N\rangle vdu\\
&\quad\quad\quad-\frac{1}{2}\int_{0}^{2\pi}\langle\frac{\partial\gamma}{\partial u} ,G(k)k\xi\rangle du
-\frac{1}{2}\int_{0}^{2\pi}\langle\gamma ,G(k)k^{2} N\rangle vdu\\
&=-\frac{1}{2}\int_{0}^{2\pi}G(k)kvdu-\frac{1}{2}\int_{0}^{2\pi}\langle k\xi ,G(k)k\xi\rangle vdu\\
&=-\frac{1}{2}\int_{0}^{2\pi}G(k)kvdu-\frac{1}{2}\int_{0}^{2\pi}G(k)k vdu\\
&=-\int_{0}^{2\pi}G(k)k vdu\\
&=-\int_{0}^{L}G(k)k ds .\endaligned
\end{equation*}

vii) \ By (iii) and (v) we drive the following equation
$$\frac{\partial k}{\partial t}=\frac{\partial}{\partial t}\frac{\partial\theta}{\partial s}
=\frac{\partial}{\partial s}\frac{\partial\theta}{\partial t}+G(k)k^{2}\frac{\partial\theta}{\partial s}
=\frac{\partial}{\partial s}\left((G(k)k)'\frac{\partial k}{\partial s}\right)+G(k)k^{3}.$$
Thus, the proof of Lemma 2.1 is completed.
\end{proof}
Similarly to \cite{GH}, we can use the angle $\theta$ of the tangent line
as the parameter of the curve and then write the curvature
$k=k(\theta)$ in terms of this parameter  which is $2\pi$ periodic
curvature function  of convex curve.  The following results gives the necessary and sufficient condition
for some one-parameter function as the curvature function of a simple closed curve(cf.
Lemma 4.1.1 in \cite{GH}).
\begin{lemma}\label{l2.2}
A positive 2$\pi$ periodic function represents the curvature
function of a  closed  and strictly convex $C^{2}$ curves in the plane if and
only if
$$\int^{2\pi}_{0}\frac{\cos\theta}{k(\theta)}d\theta=\int^{2\pi}_{0}\frac{\sin\theta}{k(\theta)}d\theta=0 .$$
\end{lemma}
According to the flow (\ref{e1.3}), we take $\tau =t$ as the time parameter
and use $\theta$ as other coordinate  and hence change variables from $(u,t)$ to $(\theta,\tau)$.
\begin{lemma}\label{l2.3}
$$\aligned\frac{\partial
k}{\partial\tau}&=k^{2}\left(\frac{\partial^{2}}{\partial\theta^{2}}(G(k)k)+G(k)k\right)\\
&=k^{2} (G(k)k)'\frac{\partial^{2}k}{\partial\theta^{2}}+k^{2}
(G(k)k)''(\frac{\partial k}{\partial\theta})^{2}+G(k)k^{3}
.\endaligned$$
\end{lemma}
\begin{proof}
By the chain rule and Lemma \ref{l2.1} (iii),(v), we get
\begin{equation*}\aligned
\frac{\partial k}{\partial t}&=\frac{\partial k}{\partial \tau}+
\frac{\partial k}{\partial \theta}\frac{\partial\theta}{\partial
t}\\&=\frac{\partial k}{\partial \tau}+ \frac{\partial k}{\partial
\theta}(G(k)k)'\frac{\partial k}{\partial s}\\&=\frac{\partial
k}{\partial \tau}+ \frac{\partial k}{\partial
\theta}(G(k)k)'\frac{\partial k}{\partial \theta}\frac{\partial
\theta}{\partial s} \\&=\frac{\partial k}{\partial
\tau}+(G(k)k)'k(\frac{\partial k}{\partial \theta})^{2} .\endaligned
\end{equation*}
On the other hand, from Lemma \ref{l2.1} (vii)  we obtain
\begin{equation*}\aligned
\frac{\partial k}{\partial t}
&=\frac{\partial}{\partial s}\left((G(k)k)'\frac{\partial k}{\partial s}\right)+G(k)k^{3}\\
&=\frac{\partial\theta}{\partial s}\frac{\partial}{\partial \theta}\left((G(k)k)'\frac{\partial k}{\partial \theta}\frac{\partial\theta}{\partial s}\right)
+G(k)k^{3}\\
&=k\frac{\partial}{\partial \theta}\left((G(k)k)'\frac{\partial k}{\partial \theta}k\right)+G(k)k^{3}\\
& =k\frac{\partial}{\partial\theta}\left(k\frac{\partial }{\partial \theta}(G(k)k)\right)+G(k)k^{3}\\
&=k^{2}\frac{\partial^{2}}{\partial\theta^{2}}(G(k)k)+(G(k)k)'k(\frac{\partial k}{\partial\theta})^{2}+G(k)k^{3}.
\endaligned
\end{equation*}
By comparing the above two equalities  we have the desired results.
\end{proof}
Throughout this paper,  we will deal with this equation and replace
$\tau$ by $t$.
\begin{lemma}\label{l2.4}
Suppose ($H1$) hold. Then the general curve shortening problem
(\ref{e1.3}) for convex curves is equivalent to the cauchy problem
\begin{equation}\label{e2.2}
\left\{ \begin{aligned} \frac{\partial k}{\partial
t}&=k^{2}\left(\frac{\partial^{2}}{\partial\theta^{2}}(G(k)k)+G(k)k\right),\quad
&\theta \in \mathbb{S}^{1},\quad &t\in[0, T),\\
k(\theta,0)&=k_{0}(\theta),
&\theta\in\mathbb{S}^{1},\quad &t=0.
\end{aligned} \right.
\end{equation}
where $0<\alpha<1$, $k\in
C^{2+\alpha,1+\frac{\alpha}{2}}(\mathbb{S}^{1}\times (0,T))$,
$k_{0}(\theta)$ is the curvature function of the initial curve
$\gamma_{0}(\theta)$.
\end{lemma}
\begin{proof}
If $\gamma(\cdot,t)$ are the classical solution of problem (\ref{e1.3}). Then by Lemma \ref{l2.3}
 the curvature function, expressed in $\theta$ coordinates, satisfies (\ref{e2.2}).

If $k_{0}(\theta)$ is the curvature function for the curve $\gamma_{0}(\theta)$ and $k(\theta,t)$ satisfies
(\ref{e2.2}). Then for each $t\geq 0$, we can define the curves by the formula
\begin{equation}\label{e2.3}
x(\theta,t)=\int^{\theta}_{0}\frac{\cos u}{k(u,t)}du,\quad
y(\theta,t)=\int^{\theta}_{0}\frac{\sin u}{k(u,t)}du.
\end{equation}
 Let $\gamma(\theta,t)=({x(\theta,t)}, {y(\theta,t)})$ , then
 $\zeta=(\cos\theta,\sin\theta)$ and $N=(-\sin\theta,\cos\theta)$ are  respectively the tangent vector and
 the inward normal vector of
the curve $\gamma(\cdot,t)$. Combining (\ref{e2.2}) with (\ref{e2.3}) we have
\begin{equation*}\aligned
\frac{\partial x}{\partial t}
&=-\int^{\theta}_{0}\frac{\cos u}{k^{2}}\frac{\partial k}{\partial t}d u \\
&=-\int^{\theta}_{0}\cos u \ \left(\frac{\partial^{2}}{\partial u^{2}}(G(k)k)+G(k)k\right)d u\\
&=-\int^{\theta}_{0}\cos u \ \frac{\partial^{2}}{\partial u^{2}}(G(k)k)d u-\int^{\theta}_{0}\cos u G(k)kd u.\qquad\quad \endaligned
\end{equation*}
So that
$$\aligned \frac{\partial x}{\partial t}
&=\!-\int^{\theta}_{0}\sin u\frac{\partial}{\partial u}(G(k)k)d u-\int^{\theta}_{0}\cos u
G(k)kd u-\cos\theta\frac{\partial}{\partial\theta}(G(k)k)+\frac{\partial}{\partial\theta}(G(k)k)|_{\theta=0}\\
&=\!\int^{\theta}_{0}\!\cos u G(k)kd u\!-\!G(k)k\sin\theta\!-\!\int^{\theta}_{0}\!\cos u
G(k)kd u\!-\!\cos\theta\frac{\partial}{\partial\theta}(G(k)k)\!+\!\frac{\partial}{\partial\theta}(G(k)k)|_{\theta=0}\\
&=-G(k)k\sin\theta-\cos\theta\frac{\partial}{\partial\theta}(G(k)k)+\frac{\partial}{\partial\theta}(G(k)k)|_{\theta=0}
.\endaligned$$
For the same reason the following equality holds
$$\frac{\partial y}{\partial t}=G(k)k\cos\theta-\sin\theta\frac{\partial}{\partial\theta}(G(k)k)-G(k)k|_{\theta=0}.$$
By setting $\theta=2\pi$ to the above two equalities we obtain
$\displaystyle\frac{\partial x}{\partial t}= 0$,
$\displaystyle\frac{\partial y}{\partial t}= 0$. On the other hand,
$x(2\pi,0)=0$, $y(2\pi,0)=0$. Because $\gamma_{0}$ is the closed
curve in the plane, then applying Lemma \ref{e2.2} we obtain $x(2\pi,0)=0$, $y(2\pi,0)=0$, and
\begin{equation*}\label{1.03}
\int^{2\pi}_{0}\frac{\cos
\theta}{k(\theta,t)}d\theta=\int^{2\pi}_{0}\frac{\sin
\theta}{k(\theta,t)}d\theta=0 .
\end{equation*}
By applying Lemma \ref{e2.2} again, the curve $\gamma(\theta,t)$ which is
defined by (\ref{e2.3}) is closed  and then we have
\begin{equation}\label{e2.4}
\frac{\partial \gamma}{\partial t}(\theta,
t)=G(k)kN-\frac{\partial}{\partial\theta}(G(k)k)\xi-(a(0,t),b(0,t)),
\end{equation}
where
$$(a(0,t),b(0,t))=\displaystyle(-\frac{\partial}{\partial\theta}(G(k)k)|_{\theta=0},G(k)k|_{\theta=0}).$$
Set
$$\theta=\theta(u,\tau),\quad t=\tau,$$
\begin{equation}\label{e2.5}
\hat{\gamma}(u,\tau)=\gamma(\theta(u,\tau),\tau)+(\int^{\tau}_{0}a(0,t)dt,\int^{\tau}_{0}b(0,t)dt),
\end{equation}
where $\theta=\theta(u,\tau)$ is  the unique  solution of the
following ordinary equation
\begin{equation*}\label{1.03}
\left\{ \begin{aligned} &\frac{\partial \theta}{\partial
\tau}=k\frac{\partial}{\partial\theta}(G(k)k) , \quad \tau\in [0,T) ,\\
&\theta(u,0)=u  .
\end{aligned} \right.
\end{equation*}
Combining  (\ref{e2.4}) with(\ref{e2.5})  we know that $\hat{\gamma}(u,\tau)$
satisfies
$$\frac{\partial\hat{\gamma}}{\partial\tau}=G(k)kN  $$ and hence we
obtain the general curve shortening flow (\ref{e1.3}).
\end{proof}
Using  the standard results on parabolic equations (cf. \cite{Kr}) we
obtain the existence results  of the problem (\ref{e2.2}).
\begin{lemma}\label{l2.5}
Suppose $k_{0}(\theta)\in C(\mathbb{S}^{1})$ and $(H1)$ holds. Then there exists
$$T>0, \quad k\in C^{2,1}(\mathbb{S}^{1}\times (0,T))\cap C(\mathbb{S}^{1}\times[0,T)),$$
satisfying (\ref{e2.2}).
\end{lemma}
By maximum principle it shows that if the initial curve are strictly
convex, then the curves remain so under the flow (\ref{e1.3}).
\begin{lemma}\label{l2.6}
Suppose ($H1$) holds and $k(\theta,t)$ is the classical solution
of (\ref{e2.2}). If $k_{0}(\theta)$ is positive on $\mathbb{S}^{1}$, then
$k_{\min}(t)=\inf \{k(\theta,t)\mid \theta\in \mathbb{S}^{1} \}$ is
a nondecreasing function.
\end{lemma}
\begin{proof}
By contradiction, there exist $0\leq t'_{1}\leq t_{1}<T$ such that
$k_{\min}(t)$ is nondecreasing in $[0,t'_{1}]$ and
$k_{\min}(t'_{1})>k_{\min}(t_{1})>0$. We suppose $t'_{1}=0$ without
losing the generalization. Set
$$t_{0}=\inf \{t\in [0,t_{1}]|k_{\min}(t)\leq k_{\min}(t_{1})\} .$$
By the regularity of $k(\theta,t)$ we know that there is
$\theta_{0}\in \mathbb{S}^{1}$ such that
$k_{\min}(t_{0})=k(\theta_{0},t_{0})$. It is easy to see that
$k_{\min}(t_{0})>0$ and then we have
$$\frac{\partial^{2}k}{\partial\theta^{2}}(\theta_{0},t_{0})\geq
0,\quad \frac{\partial k}{\partial\theta}(\theta_{0},t_{0})=0,\quad
k(\theta_{0},t_{0})>0 .$$ Hence  from (\ref{e2.2}) and ($H1$) this
yields
$$\frac{\partial k}{\partial t}(\theta_{0},t)|_{t=t_{0}}>0 ,$$
and it contradicts to the hypothesis of $t_{0}$. So that  we obtain
the desired results.
\end{proof}
Some further consequences of  Lemmas \ref{l2.4}-\ref{l2.6} are part of Proposition \ref{p1.1}.
\begin{corollary}\label{de 2.07}
 Suppose ($H1$) holds. Then there exist $T>0$ and the unique $\gamma(u,t) \in
C^{2,1}(\mathbb{S}^{1}\times (0,T))\cap C(\mathbb{S}^{1}\times
[0,T)) $ satisfying the generalized curve shortening flow (\ref{e1.3}).
\end{corollary}\label{de 2.07}
According to the flow (\ref{e1.3}) we consider the support function (cf.
\cite{S}, \cite{U}) of
 $\gamma(u,t)$ by defining
\begin{equation*}\label{1.03}
h(\theta,t)=\langle\gamma(u(\theta,t),t),(\sin\theta,-\cos\theta)\rangle,
\quad \theta\in [0,2\pi]  , \end{equation*}\label{1.03} where
$\theta$ is the tangent angle of $\gamma(\cdot,t)$ and the unit
normal vector $N=-(\sin\theta,-\cos\theta)$.

Applying the equation (\ref{e1.3}) we see that
\begin{equation*}\label{1.03}
\frac{\partial h}{\partial t}=\langle\frac{\partial \gamma}{\partial
t}+\frac{\partial u}{\partial t}\frac{\partial \gamma}{\partial
u},-N\rangle=\langle G(k)kN+\frac{\partial u}{\partial
t}|\frac{\partial \gamma}{\partial u}|\xi,-N\rangle=-G(k)k .
\end{equation*}
Using the methods in \cite{U},  we can compute the inward curvature of
$\gamma(\cdot,t)$ by the formula
\begin{equation*}
k=\left(\frac{\partial^{2}h}{\partial\theta^{2}}+h\right)^{-1}.
\end{equation*}
Then $h(\theta,t)$ satisfies the following equation
\begin{equation}\label{e2.6} \frac{\partial h}{\partial
t}\left(\frac{\partial^{2}h}{\partial\theta^{2}}
+h\right)G^{-1}\left(\left(\frac{\partial^{2}h}{\partial\theta^{2}}+h\right)^{-1}\right)=-1.
\end{equation}
By making use of the maximum principle, we can obtain the
containment principle of the flow (\ref{e1.3}) (cf. [Z]).
\begin{lemma}\label{l2.8}
 Let $\gamma_{1}$ and $\gamma_{2}$:
$\mathbb{S}^{1}\times [0,T)$ be two classical solutions of the flow
(\ref{e1.3}). If $\gamma_{2}(\cdot,0)$ is in the domain enclosed by
$\gamma_{1}(\cdot,0)$, then $\gamma_{2}(\cdot,t)$ is contained in
the domain enclosed by $\gamma_{2}(\cdot,t)$ for all $t\in [0,T)$.
\end{lemma}
\begin{proof}
 Set $h_{1}(\theta,t)$ and $h_{2}(\theta,t)$ be the support functions of $\gamma_{1}(\cdot,t)$ and
 $\gamma_{2}(\cdot,t)$. Then $h_{1}(\theta,t)$ and $h_{2}(\theta,t)$ satisfy the equations (\ref{e2.6}).
 Because $\gamma_{2}(\cdot,0)$ is in the domain enclosed by
$\gamma_{1}(\cdot,0)$  we can select $h_{1}(\theta,0)$ and $h_{2
}(\theta,0)$ such that $h_{1}(\theta,0)\geq h_{2}(\theta,0) $  for
$\theta\in \mathbb{S}^{1}$.
 Thus by applying the maximum principle of parabolic equations, we
 deduce that $h_{1}(\theta,t)\geq h_{2}(\theta,t)$ for all $t\in
 [0,T)$ and then we obtain the desired results.
\end{proof}
In order to prove some isometric inequalities we need the following
lemma which belongs to Ben.Andrews (cf, Lemma $\mathbb{I}$ 3.3 in
\cite{B}).
\begin{lemma}\label{l2.9}
Let $M$ be a compact manifold with a volume form $d\mu$, and let
$\xi$ be a continuous function on $M$. Then for any non-decreasing
function $F$, there holds
\begin{equation*}
\frac{\int_{M}\xi F(\xi)d\mu}{\int_{M} F(\xi)d\mu}\geq
\frac{\int_{M} \xi d\mu}{\int_{M} d\mu}.
\end{equation*}
\end{lemma}
The next two lemmas roughly characterized the behavior of the
geometric quantity when   $\gamma(\cdot,t)$  is contracting  to a
point under the flow (\ref{e1.3}).
\begin{lemma}\label{l2.10}
Suppose ($H1$) hold. Then under  the flow (\ref{e1.3}) we have
$$\frac{d}{d t}(\frac{L^{2}}{A})\leq 0. $$
\end{lemma}
\begin{proof}
By Lemma \ref{l2.1} (ii) and (vi), i.e,
$$\frac{dL}{dt}=-\int^{L}_{0}G(k)k^{2}ds=-\int^{2\pi}_{0}G(k)kd\theta,$$
$$\frac{dA}{dt}=-\int^{L}_{0}G(k)kds=-\int^{2\pi}_{0}G(k)d\theta, $$
we obtain
\begin{equation}\label{e2.7}
\frac{d}{dt}(\frac{L^{2}}{A})=-\frac{2L}{A}\left(\int^{2\pi}_{0}G(k)kd\theta-\frac{L}{2A}\int^{2\pi}_{0}G(k)d\theta\right).
\end{equation}
From the isometric inequality in [O]  the following inequality holds
for convex curves,
\begin{equation}\label{e2.8}
\frac{\pi L}{A}\leq \int^{L}_{0}k^{2}ds=\int^{2\pi}_{0}kd\theta,
\end{equation}
Substituting (\ref{e2.8}) into (\ref{e2.7}) we have
\begin{equation}\aligned\label{e2.9}
\frac{d}{dt}(\frac{L^{2}}{A})
&\leq -\frac{2L}{A}\displaystyle\left(\int^{2\pi}_{0}G(k)kd\theta-\frac1{2\pi}\int^{2\pi}_{0}G(k)d\theta\int^{2\pi}_{0}kd\theta\right) \\
&=-\frac{2L}{A}\int^{2\pi}_{0}G(k)d\theta\left(\displaystyle\frac{\int^{2\pi}_{0}G(k)kd\theta}{\int^{2\pi}_{0}G(k)d\theta}-
\frac{\int^{2\pi}_{0}kd\theta}{\int^{2\pi}_{0}1d\theta}\right).
\endaligned\end{equation}
Setting  $M=\mathbb{S}^{1}$, $\psi=k, d\mu=d\theta$ in Lemma 2.9
one can show that
\begin{equation}\label{e2.10}
\displaystyle\frac{\int^{2\pi}_{0}G(k)kd\theta}{\int^{2\pi}_{0}G(k)d\theta}\geq
\frac{\int^{2\pi}_{0}kd\theta}{\int^{2\pi}_{0}d\theta}.
\end{equation}
So that the proof is completed by means of (\ref{e2.9}) and (\ref{e2.10}).
\end{proof}
\begin{lemma}\label{l2.11}
Suppose ($H1$) holds. Under the flow (\ref{e1.3}) if
$$\lim_{t\rightarrow \omega}A(t)=0,$$
we have
\begin{equation}\label{e2.11}
\liminf_{t\rightarrow \omega} L\left(\int^{L}_{0}k^{2}ds-\frac{\pi L}{A}\right)\leq 0 .
\end{equation}
\end{lemma}
\begin{proof}
From (\ref{e2.7}) we see that
\begin{equation}\label{e2.12}
\aligned\frac{d}{d t}(\frac{L^{2}}{A})
&=-\frac{2L}{A}\left(\int^{2\pi}_{0}G(k)kd\theta-\frac{L}{2A}\int^{2\pi}_{0}G(k)d\theta\right)\\
&=-\frac{\int^{2\pi}_{0}G(k)d\theta}{\pi A}\cdot
L\left(\frac{2\pi\int^{2\pi}_{0}G(k)kd\theta}{\int^{2\pi}_{0}G(k)d\theta}-
\frac{\pi L}{A}\right).\endaligned
\end{equation}
Applying
$$\frac{d A}{d
t}=-\int^{2\pi}_{0}G(k)d\theta $$ to (\ref{e2.12}) we conclude that
\begin{equation*} \frac{d}{d
t}(\frac{L^{2}}{A})=\frac{1}{\pi}\frac{d}{dt}(\ln A) \cdot
L\left(\frac{2\pi\int^{2\pi}_{0}G(k)kd\theta}{\int^{2\pi}_{0}G(k)d\theta}-
\frac{\pi L}{A}\right).
\end{equation*}
Using (\ref{e2.10})  we have
\begin{equation}\label{e2.13}
\frac{d}{d t}(\frac{L^{2}}{A})\leq\frac{1}{\pi}\frac{d}{dt}(\ln
A)\cdot L\left(\int^{2\pi}_{0}kd\theta- \frac{\pi L}{A}\right).
\end{equation}
Now we prove (\ref{e2.11}) by contradiction. If not, there exists
$\delta>0$ such that if $\lim_{t\rightarrow \omega}A(t)=0$ then
\begin{equation*}
\liminf_{t\rightarrow \omega} L\left(\int^{L}_{0}k^{2}ds-\frac{\pi L}{A}\right)\geq2\delta.
\end{equation*}
Hence there exists $\beta=\beta (\delta)\in (0,\omega)$ such that if
$t\in (\beta,\omega)$, then the following inequality holds
\begin{equation}\label{e2.14}
 L\left(\int^{L}_{0}k^{2}ds-\frac{\pi L}{A}\right)\geq\delta.
\end{equation}
From (\ref{e2.13}) and (\ref{e2.14}) we obtain
\begin{equation*}
\frac{d}{d t}(\frac{L^{2}}{A})\leq\frac{\delta}{\pi}\frac{d}{dt}(\ln
A),\quad t\in (\beta,\omega).
\end{equation*}
Integrating from $\beta$ to $t$ we have
\begin{equation*}
\frac{L^{2}}{A}(t)-\frac{L^{2}}{A}(\beta)\leq\frac{\delta}{\pi}\left(\ln A(t)-\ln A(\beta)\right),
\end{equation*}
\begin{equation*}
-\frac{L^{2}}{A}(\beta)\leq\frac{\delta}{\pi}\left(\ln A(t)-\ln A(\beta)\right).
\end{equation*}
Using $\lim_{t\rightarrow \omega}A(t)=0$ one easily verifies that
\begin{equation*}
-\frac{L^{2}(\beta)}{A(\beta)}+\frac{\delta}{\pi}\ln A(\beta)=-\infty
\end{equation*}
and this contradicts to Proposition \ref{p1.1} and then the proof of Lemma \ref{l2.11} is completed.
\end{proof}

\section{asymptotic behavior}

In this section we will study the asymptotic behavior of the
curvature under the flow (\ref{e1.3}) and prove the main theorem of this paper.

We recall the following  two auxiliary results, which belong to
M.E.Gage \cite{G1}, \cite{G2} and R.Osserman \cite{O}.
\begin{lemma}(M.E.Gage)\label{l3.1}\\
 (a) There is a non-negative functional $F(\gamma)$ which is defined for
all $C^{2}$ convex curves  and  satisfies
\begin{equation}\label{e3.1}
(1-F(\gamma))\int^{L}_{0}k^{2}ds-\frac{\pi L}{A}\geq 0.
\end{equation}
(b) Given a sequence of convex curves $\{\gamma_{i}\}$ such that
$\lim_{i\rightarrow\infty}F(\gamma_{i})=0$.  If these normalized
curves $\displaystyle\eta_{i}=\sqrt{\frac{\pi}{A}}\gamma_{i}$ lie in a fixed
bounded region of the plane, then the domain $H_{i}$ which  enclosed
by $\eta_{i}$  converges to the disk in the Hausdorff metric. \\
(c) $F(\gamma)$=0 if and only if $\gamma$ is a circle.
\end{lemma}

\begin{lemma}(Bonneson inequality)\label{l3.2} Let $\gamma$
be a $C^{1}$ closed convex curve. Then
\begin{equation}\label{e3.2}
\frac{L^{2}}{A}-4\pi\geq \frac{\pi^{2}}{A}(r_{out}-r_{in})^{2} .
\end{equation}
\end{lemma}

\begin{definition}\label{d3.3}
Let A, B be two closed convex sets and $A_{\epsilon}$= $\{$$x\in
{\mathbb R}^{2}\mid dist (x,A)\leq\epsilon$$\}$. Then the Hausdorff distance
between the sets A and B is given by
\begin{equation*}
d_{H}(A,B)=\inf \{\epsilon\mid A\subseteq B_{\epsilon}, B\subseteq
A_{\epsilon}\}.
\end{equation*}
\end{definition}

{\bf Proof of  Theorem 1.2 (i).} \  We use the idea in \cite{G2}.  Consider the geometric quantities $L(t)$ and $A(t)$ according to the
flow $\gamma(\cdot,t)$ satisfying (\ref{e1.3}). By (\ref{e3.1}) we have
\begin{equation}\label{e3.3}
\int^{L(t)}_{0}k^{2}(\theta,t)ds-\frac{\pi L(t)}{A(t)}\geq
F(\gamma(t))\int^{L(t)}_{0}k^{2}(\theta,t)ds.
\end{equation}
Using Schwartz inequality we see that
\begin{equation}\label{e3.4}
(2\pi)^{2}=\left(\int^{L(t)}_{0}k(\theta,t)ds\right)^{2}\leq
\int^{L(t)}_{0}k^{2}(\theta,t)ds\int^{L(t)}_{0}1ds=L(t)\int^{L(t)}_{0}k^{2}(\theta,t)ds.
\end{equation}
Substituting (\ref{e3.4}) into (\ref{e3.3}) we obtain
\begin{equation}\label{e3.5}
L(t)\left(\int^{L(t)}_{0}k^{2}(\theta,t)ds-\frac{\pi L(t)}{A(t)}\right)\geq
F(\gamma)L(t)\int^{L(t)}_{0}k^{2}(\theta,t)ds\geq
4\pi^{2}F(\gamma(t)).
\end{equation}

In the following steps we will show that if $\forall t_{i}\in
[0,\omega)$ which satisfies $\lim_{i\rightarrow
+\infty}t_{i}=\omega$ then the limitations in Theorem 1.2 (i),(ii) hold.

By substituting
$\gamma_{i}=\gamma(\cdot,t_{i})$, $A(t_{i})$, $L(t_{i})$ into (\ref{e3.5})
and using (\ref{e2.11})  we have
\begin{equation*}
\lim_{i\rightarrow+\infty}F(\gamma_{i})=0.
\end{equation*}
Next we show that the normalized curve
$\displaystyle\eta_{i}=\sqrt{\frac{\pi}{A}}\gamma_{i}$ lies in a bound region.
From Lemma \ref{l2.10}, we observe that $\displaystyle\frac{L^{2}}{A}$ decreases under
under the flow (\ref{e1.3}). One easily verifies that
$$\frac{L^{2}(t_{i})}{A(t_{i})}=\frac{L^{2}_{n}(t_{i})}{A_{n}(t_{i})},$$
and
\begin{equation}\label{e3.6}
\frac{L^{2}(t_i)}{A(t_i)}-4\pi=\frac{L^{2}_{n}(t_i)}{A_n(t_i)}-4\pi\geq
\frac{\pi^{2}}{A_n(t_i)}(r_{n,out}(t_i)-r_{n,in}(t_i))^{2}=\pi(r_{n,out}(t_i)-r_{n,in}(t_i))^{2},
\end{equation}
and $r_{n,in}(t_i)\leq 1$, by using the Bonneson inequality, where $A_n(t_i)$ is the area about the bound
domain enclosed by the normalized curve $\eta_{i}$, and $L_{n}(t_i)$ is  the
perimeter of the curve  $\eta_{i}$, $r_{n,out}(t_i)$ and $r_{n,in}(t_i)$ are
respectively the radii of the largest and smallest circumscribed
circle of the curve $\eta_{i}$. By (\ref{e3.6}) it shows that the outer radii of
the normalized curve  $\eta_{i}$ are bounded for all $t_{i}\in
[0,\omega)$.  From Proposition \ref{p1.1} we know that $\gamma_{i}$
shrinks to a point under the flow (\ref{e1.3}). Hence  if we use
$\vartheta$ as the origin in the homothetic expansion of
$\mathbb{R}^{2}$, then all of the normalized curve  $\eta_{i}$  lies
in a ball of radius 2C around this point.

Applying Lemma \ref{l3.1}(b), we see that the sequence
of normalized domain $H(t_{i})$ according to $\eta_{i}$ converges to
the unit disk in the Haudorff metric,
\begin{equation}\label{e3.7}
\lim_{i\rightarrow+\infty}H(t_{i})=H_{0},
\end{equation}
where $H_{0}$ is the unit disk in the plane.

Since L and A are continuous functionals of convex domain then there holds,
\begin{equation}\label{e3.8}
\lim_{i\rightarrow+\infty}\frac{L^{2}(t_{i})}{A(t_{i})}=\lim_{i\rightarrow+\infty}\frac{L_{n}^{2}(t_{i})}{A_{n}(t_{i})}=
\lim_{i\rightarrow+\infty}\frac{L^{2}(H(t_{i}))}{A(H(t_{i}))}=
\frac{L^{2}(H_{0})}{A(H_{0})}=4\pi .\end{equation}
From (\ref{e3.2}) we have
\begin{equation*}
\frac{L^{2}(t_{i})}{A(t_{i})}-4\pi\geq
\frac{\pi^{2}}{A(t_{i})}\left(r_{out}(t_{i})-r_{in}(t_{i})\right)^{2}\geq
\frac{\pi^{2}r^{2}_{out}(t_{i})}{A(t_{i})}\left(1-\frac{r_{in}(t_{i})}{r_{out}(t_{i})}\right)^{2}.
\end{equation*}
It is easy to see that $\pi r^{2}_{out}(t_{i})\geq A(t_{i})$, and then
\begin{equation*}
\frac{L^{2}(t_{i})}{A(t_{i})}-4\pi\geq
\left(1-\frac{r_{in}(t_{i})}{r_{out}(t_{i})}\right)^{2}.
\end{equation*}
Such that combining this with (\ref{e3.8}) we have
\begin{equation}\label{e3.9}
\lim_{i\rightarrow +\infty}\frac{r_{in}(t_{i})}{r_{out}(t_{i})}=1.
\end{equation}
 \hfill$\Box$

To prove  Theorem 2.1 (ii),(iii), we need the following gradient
estimates of the curvature. A similar proof can be found in \cite{Z}.

\begin{lemma}\label{l3.4}
Set $\Phi(k)=G(k)k$ and let $k=k(\theta,t)$ be the curvature function of
the flow (\ref{e1.3}), where $\theta$ is the tangent angle of the curve
$\gamma(\cdot,t)$. Suppose ($H1$) , ($H2$) hold and $\varpi\in
(0,\omega)$. Then the following inequality holds,
\begin{equation}\label{e3.10}
\max_{0\leq t\leq \varpi,\theta\in \mathbb{S}^{1}}\mid\frac{\partial\Phi}{\partial\theta}\mid^{2}
\leq\max\left\{2\max_{0\leq t\leq\varpi,\theta\in \mathbb{S}^{1} }\Phi^{2},
\max_{t=0,\theta\in \mathbb{S}^{1}}\left(\mid\frac{\partial \Phi}{\partial\theta}\mid^{2}+2\Phi^{2}\right)\right\} .
\end{equation}
\end{lemma}
\begin{proof}
By Lemma \ref{l2.6} we know that $k_{\min}(t)>0$ for $t\in [0,\omega)$. It
follows from Lemma \ref{l2.3} that $\phi=\Phi(k)$ satisfies the following
equation
\begin{equation}\label{e3.11}
\frac{\partial \Phi}{\partial
t}=k^{2}\Phi'\frac{\partial^{2}\Phi}{\partial
\theta^{2}}+k^{2}\Phi'\Phi.
\end{equation}
Set
\begin{equation*} \Psi=(\frac{\partial
\Phi}{\partial\theta})^{2}+\lambda \Phi^{2},
\end{equation*}
where $\lambda$ is a constant to be determined. Suppose
$(\theta_{0}, t_{0})\in \mathbb{S}^{1}\times(0,\varpi]$ such that
\begin{equation*}
\Psi(\theta_{0},t_{0})=\max_{\mathbb{S}^{1}\times
[0,\varpi]}\left((\frac{\partial \Phi}{\partial\theta})^{2}+\lambda
\Phi^{2}\right).
\end{equation*}
Then at $(\theta_{0}, t_{0})$, $\Psi$ satisfies the following
properties
\begin{equation}\label{e3.12}
\frac{\partial\Psi}{\partial\theta}=0,\quad
\frac{\partial^{2}\Psi}{\partial\theta^{2}}\leq 0,\quad
\frac{\partial\Psi}{\partial t}\geq 0.
\end{equation}
Next we will prove that if  selecting  some constant $\lambda$ so
large, then at $(\theta_{0}, t_{0})$ we have
\begin{equation}\label{e3.13}
\frac{\partial\Phi}{\partial\theta}=0 .
\end{equation}
Suppose not, then using
\begin{equation*}
0=\frac{\partial \Psi}{\partial\theta}=2\frac{\partial
\Phi}{\partial\theta}\left(\frac{\partial^{2}
\Phi}{\partial\theta^{2}}+\lambda \Phi\right)
\end{equation*}
we see that
\begin{equation}\label{e3.14}
0=\frac{\partial^{2} \Phi}{\partial\theta^{2}}+\lambda \Phi.
\end{equation}
From (\ref{e3.11}) and (\ref{e3.12}) we have
\begin{equation}\label{e3.15}
\aligned 0& \leq\frac{1}{2}\frac{\partial \Psi}{\partial t} \\
&=\frac{\partial\Phi}{\partial\theta}\frac{\partial^{2}\Phi}{\partial \theta\partial
t}+\lambda \Phi\frac{\partial \Phi}{\partial t}\\
&=\frac{\partial}{\partial\theta}(\Phi'k^{2})\frac{\partial
\Phi}{\partial\theta}\frac{\partial^{2}
\Phi}{\partial\theta^{2}}+\Phi'k^{2}\frac{\partial
\Phi}{\partial\theta}\frac{\partial^{3}
\Phi}{\partial\theta^{3}}+\frac{\partial}{\partial\theta}(\Phi'k^{2})
\Phi\frac{\partial \Phi}{\partial\theta }\\
& \ \ \ \ +k^{2}\Phi'(k)(\frac{\partial \Phi}{\partial \theta})^{2}+\lambda
\Phi\Phi'k^{2}\frac{\partial^{2} \Phi}{\partial\theta^{2}}+\lambda
k^{2}\Phi'\Phi^{2}.
\endaligned\end{equation}
By $\displaystyle\frac{\partial^{2}\Psi}{\partial\theta^{2}}\leq 0$
 we see that
\begin{equation}\label{e3.16}
0\geq \frac{\partial
\Phi}{\partial\theta}\left(\frac{\partial^{3}\Phi}{\partial\theta^{3}}+
\lambda\frac{\partial \Phi}{\partial\theta}\right).
\end{equation}
Substituting (\ref{e3.14}), (\ref{e3.16}) into (\ref{e3.15}) we obtain
\begin{equation}\label{e3.17}
\aligned
0&\leq-\lambda\frac{\partial}{\partial\theta}(\Phi'k^{2})\Phi\frac{\partial \Phi}{\partial\theta}-\lambda k^{2}
\Phi'(\frac{\partial\Phi}{\partial\theta})^{2}+\frac{\partial}{\partial\theta}(\Phi'k^{2})\Phi\frac{\partial \Phi}{\partial\theta} \\
&\quad\quad\quad +k^{2}\Phi'(\frac{\partial\Phi}{\partial \theta})^{2}-\lambda^{2} k^{2}\Phi^{2}\Phi'+ \lambda k^{2} \Phi'\Phi^{2} \\
&=(1-\lambda)\frac{\partial}{\partial\theta}(\Phi'k^{2})\Phi\frac{\partial \Phi}{\partial\theta}+(1-\lambda)\Phi'k^{2}
(\frac{\partial\Phi}{\partial\theta})^{2}+(\lambda-\lambda^{2})k^{2} \Phi'\Phi^{2}.
\endaligned\end{equation}
By the definition of $\Phi$ and ($H1$), ($H2$), we have
\begin{equation*}
\Phi>0,\qquad \Phi'>0,
\end{equation*}
\begin{equation*}
(\Phi'k^{2})'=G''(k)k^{3}+4k^{2}G'(k)+2G(k)k=k(G(k)k^{2})''\geq 0,
\end{equation*}
\begin{equation*}
\frac{\partial}{\partial\theta}(\Phi'k^{2})\Phi\frac{\partial
\Phi}{\partial\theta}= (\Phi'k^{2})'(\frac{\partial k}{\partial
\theta})^{2}\Phi'\Phi\geq 0,
\end{equation*}
\begin{equation}\label{e3.18}
\Phi'k^{2} \Phi(\frac{\partial \Phi}{\partial\theta})^{2}>0,\quad
k^{2} \Phi'\Phi^{2}>0.
\end{equation}
 By selecting $\lambda=2$ and hence substituting it into (\ref{e3.17}) and using
(\ref{e3.18}), we obtain the contradiction. Such  that  (\ref{e3.13}) holds.

By (\ref{e3.13}) we arrive at
\begin{equation*}
\aligned
\max_{0\leq t\leq \varpi,\theta\in \mathbb{S}^{1}}\mid\frac{\partial \Phi}{\partial\theta}\mid^{2}
&\leq \max_{\mathbb{S}^{1}\times [0,\varpi]}\left((\frac{\partial \Phi}{\partial\theta})^{2}+\lambda\Phi^{2}\right) \\
&=\Psi(\theta_{0},t_{0}) \\
&=2\Phi^{2}|_{\theta=\theta_{0},t=t_{0}} \\
&\leq\max\left\{2\max_{0\leq t\leq \varpi,\theta\in\mathbb{S}^{1} }\Phi^{2},
\max_{t=0,\theta\in\mathbb{S}^{1}}\left(\mid\frac{\partial\Phi}{\partial\theta}\mid^{2}+2\Phi^{2}\right)\right\}.
\endaligned\end{equation*}
So that  the proof is completed.
\end{proof}

\begin{lemma}\label{l3.5}
 Let $q(t)$ be continuous function on $[0,\omega)$. Suppose for each $\varpi\in
 [0,\omega)$, satisfying
 $$\sup_{0\leq t\leq\varpi}q(t)<+\infty, \ \ \ \
\lim_{t\rightarrow\omega}q(t)=+\infty.$$
Then there exists $\{t_{i}\}\subset[0,\omega)$, satisfying
\begin{equation}\label{e3.19}
 \forall i\in \{1,2,\cdots\}, t_{i}<t_{i+1},
\lim_{i\rightarrow +\infty}t_{i}=\omega,
\end{equation}
\begin{equation*}
q(t_{i})=\sup_{0\leq t\leq t_{i}}q(t) .
\end{equation*}
\end{lemma}

\begin{proof}
Consider the sequence $\{t'_{i}\}\triangleq
\{T-\displaystyle\frac{T}{i+1}\}$. Firstly select $t_{1}\in [0,
t'_{1}]$  satisfying
$$q(t_{1})=\sup_{0\leq t\leq t'_{1}}q(t).$$
Then
$$q(t_{1})=\sup_{0\leq t\leq t_{1}}q(t).$$
It follows from $\displaystyle\lim_{t\rightarrow
\omega}q(t)=+\infty$ that we can choose $t''_{j_{1}}\in \{t'_{i}\}$
satisfying $q(t''_{j_{1}})>q(t_{1})+2$. So that we can take
$t_{2}\in [0,t''_{j_{1}}]$ satisfying
$$q(t_{2})=\sup_{0\leq t\leq t''_{j_{1}}}q(t).$$
Thus
$$q(t_{2})=\sup_{0\leq t\leq t_{2}}q(t).$$
In general we can select $t''_{j_{n}}\in \{t'_{i}\}$ satisfying
$q(t''_{j_{n}})>q(t_{n})+n+1$ and then choose $t_{n+1}\in
[0,t''_{j_{n}}]$ satisfying
$$q(t_{n+1})=\sup_{0\leq t\leq t''_{j_{n}}}q(t).$$
Then  there holds
$$q(t_{n+1})=\sup_{0\leq t\leq t_{n+1}}q(t).$$
The desired results  follows by taking trace.
\end{proof}
Set $q(t)=k_{\max}(t)$. Then by Proposition \ref{p1.1}  it is
easy to verify that $k_{\max}(t)$ satisfies the conditions of Lemma \ref{l3.5}.
\begin{lemma}\label{l3.6}
Suppose ($H1$) and ($H2$) hold. If we take the sequence
$\{t_{i}\}$ satisfying (\ref{e3.19}) such that the following holds
\begin{equation*}
\forall i\in \{1,2,\cdots\},\quad k_{\max}(t_{i})=\sup_{0\leq t\leq
t_{i}}k_{\max}(t)
\end{equation*}
and for each $i\in \{1,2,\cdots\}$ there exists $ \theta_{i0}\in
\mathrm{S}^{1}$ such that $k_{\max}(t_{i})=k(\theta_{i0},t_{i})$.
Then there exist constants $C, C_{1}>0$ depending only on
$\gamma_{0}$ such that
\begin{equation}\label{e3.20}
(1-2\mid\theta-\theta_{i0}\mid)\Phi(k_{\max}(t_{i}))\leq
\Phi(k(\theta,t_{i}))+C,\quad \forall \theta\in \mathbb{S}^{1},
\end{equation}
\begin{equation}\label{e3.21}
\Phi(k_{\max}(t_{i}))\leq C_{1}\Phi(k(\theta,t_{i})), \quad
\forall\theta\in \mathbb{S}^{1},
\end{equation}
where $\Phi(x)=G(x)x.$
\end{lemma}
\begin{proof}
Step 1. Set $\Phi(\theta)=\Phi(k(\theta,t_{i})).$ For each
$\theta\in \mathbb{S}^{1}$ by the medium theorem and (\ref{e3.10}) we have
\begin{equation}\label{e3.22}
\Phi(\theta_{i0})-\Phi(\theta)=\frac{\partial\Phi}{\partial\theta}(\hat{\theta})(\theta_{i0}-\theta)
\leq (2\Phi(\theta_{i0})+C)|\theta_{i0}-\theta|\leq
2\Phi(\theta_{i0})|\theta_{i0}-\theta|+C.
\end{equation}
This yields the inequality (\ref{e3.20}).

Step 2. Take $i$ such that $\Phi(\theta_{i0})$ is large sufficiently. Being likely with (\ref{e3.22}), if $\theta_{1},
\theta\in \mathbb{S}^{1}$, we obtain
\begin{equation}\label{e3.23}
\Phi(\theta_{1})-\Phi(\theta)\leq
2\Phi(\theta_{i0})|\theta_{1}-\theta|+C.
\end{equation}
It is clear that $\Phi(\theta)=\Phi(k(\theta,t_{i}))\geq
\Phi(k_{\min}(0))>0.$ Set $\theta_{1}=\theta_{i0}$. If
$|\theta-\theta_{i0}|\leq \displaystyle\frac{1}{8}$, then from
(\ref{e3.23}) we have
\begin{equation*}
\Phi(\theta_{i0})-\Phi(\theta)\leq
2\Phi(\theta_{i0})|\theta_{i0}-\theta|+C\leq
\frac{1}{2}\Phi(\theta_{i0}).
\end{equation*}
So that
\begin{equation*}
\frac{1}{2}\Phi(\theta_{i0})\leq \Phi(\theta).
\end{equation*}
Set $\theta_{1}=\theta_{i0}+\displaystyle\frac{1}{8}$. Similarly if
$|\theta-\theta_{1}|\leq\displaystyle\frac{1}{16}$, then using
(\ref{e3.23}) and $\displaystyle\frac{1}{2}\Phi(\theta_{i0})\leq
\Phi(\theta_{i0}+\displaystyle\frac{1}{8}) $ we obtain
\begin{equation*}
\frac{1}{4}\Phi(\theta_{i0})\leq \Phi(\theta).
\end{equation*}
In general, for each $n\in\mathbb{N},$ set
$$\Lambda_{n}=
\displaystyle\biggl[\theta_{i0}+\sum^{n}_{s=1}\frac{1}{8s},\,
\theta_{i0}+\sum^{n}_{s=1}\frac{1}{8s}+\frac{1}{8(n+1)}\biggr].$$
Then if $\theta\in \Lambda_{n}$, such that there holds
\begin{equation*}
\frac{\Phi(\theta_{i0})}{2(n+1)}\leq \Phi(\theta).
\end{equation*}
Given $n$ being so large , we may cover $\mathbb{S}^{1}$ by
$\Lambda_{1},$ $\Lambda_{2},$ $\cdots,$ $\Lambda_{n}$. For
each $\theta\in \mathbb{S}^{1}$, we have
\begin{equation*}
\frac{\Phi(\theta_{i0})}{2(n+1)}\leq \Phi(\theta).
\end{equation*}
Taking $C_{1}=2(n+1)$ we have the desired results.
\end{proof}
Set
\begin{equation*} \hat{k_{\sigma}}(t_{i})=\sup\{\inf_{[a,b]}
k(\theta,t_{i})|[a,b]\subset (-\infty,+\infty), b-a=\sigma\}.
\end{equation*}
 We introduce a lemma of M.E.Gage and R.S.Hamilton (cf. Lemma 5.1
in \cite{GH}), which is crucial for studying the asymptotic behavior of
the curvature under the general curve shortening flow.

\begin{lemma}\label{l3.7}
\begin{equation*}
\hat{k}_{\sigma}(t_{i})r_{in}(t_{i})\leq
\frac{1}{1-\Lambda(\sigma)\left(\frac{r_{out}(t_{i})}{r_{in}(t_{i})}-1\right)},
\end{equation*}
where
$$\Lambda(\sigma)=\frac{2\cos\frac{\sigma}{2}}{1-\cos\frac{\sigma}{2}}.$$
\end{lemma}
\begin{remark}
The proof of Lemma 5.1 in \cite{GH} follows only from the convexity of the closed
curve $\gamma(\cdot,t_{i})$.
\end{remark}
\begin{corollary}\label{C3.9}
Suppose ($H1$) and ($H2$) hold. Consider the sequence
$\{t_{i}\}$ satisfying the conditions of Lemma \ref{l3.6}. Then for the
positive $\epsilon$ being  small sufficiently, we have
\begin{equation}\label{e3.24}
k_{\max}(t_{i})r_{in}(t_{i})\leq
\frac{2}{1-\epsilon}\cdot\frac{1}{1-\Lambda(\epsilon)\left(\frac{r_{out}(t_{i})}{r_{in}(t_{i})}-1\right)}.
\end{equation}
\end{corollary}
\begin{proof}
It follows from (\ref{e3.20}) that
\begin{equation*}
\aligned
(1-2\mid\theta-\theta_{i0}\mid)G(k_{\max}(t_{i}))k_{\max}(t_{i})&\leq
G(k(\theta,t_{i}))k(\theta,t_{i})+C\\
&\leq G(k_{\max}(t_{i}))k(\theta,t_{i})+C ,\quad \forall \theta\in
\mathbb{S}^{1}.\endaligned
\end{equation*}
So that
\begin{equation}\label{e3.25}
(1-2\mid\theta-\theta_{i0}\mid)k_{\max}(t_{i})\leq k(\theta,t_{i})+C
,\quad \forall \theta\in \mathbb{S}^{1}.
\end{equation}
By Proposition \ref{p1.1} we have
$$\lim_{i\rightarrow +\infty}k_{\max}(t_{i})=+\infty.$$
Hence from (\ref{e3.21}) and ($H1$) one can easily
verify that $$\displaystyle\lim_{i\rightarrow +\infty}k(\theta,
t_{i})=+\infty, \quad \forall \theta\in \mathbb{S}^{1}.$$ Combining
 this with (\ref{e3.25}) we obtain
\begin{equation*}
(1-2\mid\theta-\theta_{i0}\mid)k_{\max}(t_{i})\leq 2k(\theta,t_{i})
\end{equation*}
for $i$ being  large enough, and $\forall\theta\in
\mathbb{S}^{1}$. Given any  $\epsilon>0$ , if
$\mid\theta-\theta_{i0}\mid\leq \displaystyle\frac{\epsilon}{2}$
then
\begin{equation*}
2k(\theta,t_{i})\geq k_{\max}(t_{i})(1-\epsilon).
\end{equation*}
Takes $\sigma=\epsilon$. It follows from the definition of $\hat{k_{\sigma}}(t_{i})$ that we have
\begin{equation*}
2\hat{k}_{\sigma}(t_{i})\geq k_{\max}(t_{i})(1-\epsilon).
\end{equation*}
Then using Lemma \ref{l3.7} we obtain
\begin{equation*}
k_{\max}(t_{i})r_{in}(t_{i})(1-\epsilon)\leq
2\hat{k}_{\sigma}(t_{i})r_{in}(t_{i})\leq
\frac{2}{1-\Lambda(\epsilon)\left(\frac{r_{out}(t_{i})}{r_{in}(t_{i})}-1\right)}.
\end{equation*}
This yields the desired results.
\end{proof}
\begin{corollary}\label{c3.10}
Suppose ($H1$) and ($H2$) hold. Consider the sequence
$\{t_{i}\}$ satisfying the conditions of Lemma 3.6. Then for
 the positive $\epsilon$ being small sufficiently, there exists
$i(\epsilon)\in \mathbb{N}$, such that if $i>i(\epsilon)$,  we have
\begin{equation*}
k_{\max}(t_{i})r_{in}(t_{i})\leq \frac{2}{(1-\epsilon)^{2}}.
\end{equation*}
\end{corollary}

\begin{proof}
By Theorem 1.2 (ii),
\begin{equation*}
\lim_{i\rightarrow +\infty}\frac{r_{in}(t_{i})}{r_{out}(t_{i})}=1.
\end{equation*}
For $i$ being so large  we have
\begin{equation}\label{e3.26}
1-\Lambda(\epsilon)\left(\frac{r_{out}(t_{i})}{r_{in}(t_{i})}-1\right)\geq1-\epsilon.
\end{equation}
Then substituting (\ref{e3.26}) into (\ref{e3.24}) we obtain the desired results.
\end{proof}
\begin{theorem} \label{de 5.05}
Suppose ($H1$) and ($H2$) hold. Consider the sequence
$\{t_{i}\}$ satisfying the conditions of Lemma 3.6. Then we have
\begin{equation}
\lim_{i\rightarrow+\infty}k(\theta,t_{i})r_{in}(t_{i})=1,\quad
\forall \theta\in \mathbb{S}^{1}.
\end{equation}
\end{theorem}

\begin{proof}
Set $f_{i}(\theta)=k(\theta,t_{i})r_{in}(t_{i})$ and $\Phi(x)=G(x)x$.

Step 1. We will prove that $f_{i}(\theta)$ is equi-continuous and
bounded uniformly.

Because $G(x)$ is non-decreasing function for $x\in (0,+\infty)$.
Then by Lemma 3.4, for $\theta\in \mathbb{S}^{1}$, we arrive at
\begin{equation*}\aligned
|G(k(\theta,t_{i}))\frac{\partial k}{\partial\theta}(\theta,t_{i})|
&\leq|\frac{\partial\Phi}{\partial\theta}(\theta,t_{i})| \\
&\leq \max_{0\leq t\leq t_{i},\theta\in \mathbb{S}^{1}}|\frac{\partial \Phi}{\partial\theta}|\\
&\leq 2\max_{0\leq t\leq t_{i},\theta\in \mathbb{S}^{1}}|\Phi|+C \\
&=2\Phi(k_{\max}(t_{i}))+C. \endaligned
\end{equation*}
So that
\begin{equation*}
\Phi(k(\theta,t_{i}))|\frac{\partial
k}{\partial\theta}(\theta,t_{i})|=|G(k(\theta,t_{i}))k(\theta,t_{i})\frac{\partial
k}{\partial\theta}(\theta,t_{i})|\leq
2\Phi(k_{\max}(t_{i}))k_{\max}(t_{i})+Ck_{\max}(t_{i}),
\end{equation*}
\begin{equation*}
|\frac{\partial k}{\partial\theta}(\theta,t_{i})|\leq
2\frac{\Phi(k_{\max}(t_{i}))}{\Phi(k(\theta,t_{i}))}k_{\max}(t_{i})+C\frac{k_{\max}(t_{i})
}{\Phi(k(\theta,t_{i}))}\leq
2\frac{\Phi(k_{\max}(t_{i}))}{\Phi(k(\theta,t_{i}))}k_{\max}(t_{i})+C\frac{k_{\max}(t_{i})
}{\Phi(k_{\min}(0))}.
\end{equation*}
It follows from  (3.25) that
\begin{equation*} |\frac{\partial
k}{\partial\theta}(\theta,t_{i})|\leq (2C+C)k_{\max}(t_{i}) ,\quad
\forall\theta\in \mathbb{S}^{1}.
\end{equation*}
By Corollary 3.10  we obtain
\begin{equation*}
|\frac{\partial k}{\partial\theta}(\theta,t_{i})r_{in}(t_{i})|\leq
(2C+C)k_{\max}(t_{i})r_{in}(t_{i})\leq C
\end{equation*}
for $i$ being so large. This yields
\begin{equation*}
|\frac{d f_{i}}{d\theta}(\theta)|\leq C.
\end{equation*}
On the other hand, by Corollary \ref{c3.10},
\begin{equation*}
| f_{i}(\theta)|\leq C.
\end{equation*}
The proof of Step 1 is completed.

Step 2. Because $f_{i}(\theta)$ is equi-continuous and bounded
uniformly. Then by Ascoli-Arzela theorem, there exists $f(\theta)\in
C(\mathbb{S}^{1})$ such that
\begin{equation}\label{e3.28}
\lim_{i\rightarrow+\infty}f_{i}(\theta)=f(\theta),\quad
\forall\theta\in \mathbb{S}^{1}.
\end{equation}

Step 3. We will prove that $f(\theta)\leq 1$ for $\forall\theta\in
\mathbb{S}^{1}$.

Suppose the assertion is false. Then there exists $\theta_{0}\in
\mathbb{S}^{1}, \beta>0$, such that $f(\theta_{0})\geq 1+3\beta$.
Hence there exists also $\delta>0$, such that if $\theta\in
[\theta_{0}-\delta,\theta_{0}-\delta]$, we have
\begin{equation*}
f(\theta)\geq 1+2\beta.
\end{equation*}
By (\ref{e3.28}) for $i$ being so large  we have
\begin{equation*}
f_{i}(\theta)\geq 1+\beta,\quad \forall\theta\in
[\theta_{0}-\delta,\theta_{0}-\delta],
\end{equation*}
i.e,
\begin{equation*}
k(\theta,t_{i})r_{in}(t_{i})\geq 1+\beta,\quad \forall\theta\in
[\theta_{0}-\delta,\theta_{0}-\delta].
\end{equation*}
Take $\sigma=2\delta$. Then according to the definition of $\hat{k_{\sigma}}(t_{i})$ we obtain
\begin{equation*}
1+\beta\leq\hat{k}_{2\delta}(t_{i})r_{in}(t_{i})\leq
\frac{1}{1-\Lambda(2\delta)\left(\frac{r_{out}(t_{i})}{r_{in}(t_{i})}-1\right)}.
\end{equation*}
Then using \begin{equation*} \lim_{i\rightarrow
+\infty}\frac{r_{in}(t_{i})}{r_{out}(t_{i})}=1.
\end{equation*}
we have
$$1+\beta\leq 1 ,$$
and it is impossible. So $f(\theta)\leq 1$ for $\forall\theta\in \mathbb{S}^{1}$.

Step 4. We will prove that $f(\theta)\equiv 1$.

By Fatou lemma  we have
\begin{equation}\label{e3.29}
\int^{2\pi}_{0}\frac{d\theta}{f(\theta)}\leq
\liminf_{i\rightarrow+\infty}\int^{2\pi}_{0}\frac{d\theta}{f_{i}(\theta)}
=\liminf_{i\rightarrow+\infty}\int^{2\pi}_{0}\frac{d\theta}{k(\theta,t_{i})r_{in}(t_{i})}.
\end{equation}
By the convexity of $\gamma(\cdot,t_{i})$, it is easy to verify that
$$L(t_{i})=\int^{2\pi}_{0}\frac{d\theta}{k(\theta,t_{i})} $$
and substitute it into (\ref{e3.29}) we obtain
\begin{equation}\label{e3.30}
\int^{2\pi}_{0}\frac{d\theta}{f(\theta)}\leq
\liminf_{i\rightarrow+\infty}\frac{L(t_{i})}{r_{in}(t_{i})}=\liminf_{i\rightarrow+\infty}
\frac{L(t_{i})}{r_{out}(t_{i})}\cdot\frac{r_{out}(t_{i})}{r_{in}(t_{i})}.
\end{equation}
By the geometric property of $r_{out}$ one can easily verify  that
$2\pi r_{out}\geq L $. Then combining this with (\ref{e3.30}) we have
\begin{equation*}
\int^{2\pi}_{0}\frac{d\theta}{f(\theta)}\leq
\liminf_{i\rightarrow+\infty}\frac{L(t_{i})}{r_{in}(t_{i})}
\leq 2\pi\cdot\liminf_{i\rightarrow+\infty}\frac{r_{out}(t_{i})}{r_{in}(t_{i})} .
\end{equation*}
By making use of Theorem 1.2 (i) again we have
\begin{equation*}
\int^{2\pi}_{0}\frac{d\theta}{f(\theta)}\leq 2\pi.
\end{equation*}
On the other hand, by $f(\theta)\leq 1$ in Step 3 we obtain
\begin{equation*}
\int^{2\pi}_{0}\frac{d\theta}{f(\theta)}\geq 2\pi.
\end{equation*}
This yields
\begin{equation*}
\int^{2\pi}_{0}\frac{d\theta}{f(\theta)}= 2\pi.
\end{equation*}
Using $f(\theta)\leq 1$ again we have $f(\theta)\equiv 1$.

Combining (\ref{e3.28}) with Step 4 we completed the proof the theorem.
\end{proof}

\begin{remark}\label{de 2.02}
It follows from Cauchy criterion that the following limitation holds
\begin{equation}\label{e3.31}
\lim_{t\rightarrow \omega}k(\theta,t)r_{in}(t)=1,\quad
\forall\theta\in \mathbb{S}^{1}.
\end{equation}
\end{remark}
{ \bf Proof of the theorem 1.2 (ii),(iii).}

Step 1. According to (3.31) we conclude that
 \begin{equation*}
\lim_{t\rightarrow \omega}k_{\max}(t)r_{in}(t)=1,\quad
\lim_{t\rightarrow
 \omega}k_{\min}(t)r_{in}(t)=1.
\end{equation*}
Combining  this with Theorem 1.2 (ii) it shows that
\begin{equation}\label{e3.32}
\lim_{t\rightarrow \omega}\frac{k_{\max}(t)}{k_{\min}(t)}=1 .
\end{equation}

Step 2. Given $t\in (0,\omega)$ and consider $k_{\max}(t)$. By the
property of continuous function, there exists $\theta=\theta(t)\in
\mathbb{S}^{1}$ such that $k_{\max}(t)=k(\theta(t),t)$. Then  at
$(\theta(t),t)$ by the regularity of $k(\theta,t)$  we have
 \begin{equation}\label{e3.33}
  \frac{\partial k}{\partial
\theta}=0,\quad \frac{\partial^{2} k}{\partial \theta^{2}}\leq
0,\quad \frac{dk}{dt}=\frac{\partial k}{\partial
\theta}\frac{d\theta}{dt}+\frac{\partial k}{\partial
t}=\frac{\partial k}{\partial t}.
\end{equation}
By Lemma \ref{l2.3},
\begin{equation}\label{e3.34}
 \frac{\partial k}{\partial t}=k^{2}\left(\frac{\partial^{2}}{\partial\theta^{2}}(G(k)k)+G(k)k\right).
\end{equation}
Combining (\ref{e3.33}) with (\ref{e3.34}) we see that
\begin{equation*}
 \frac{dk_{\max}(t)}{d
t}\leq G(k_{\max}(t))k^{3}_{\max}(t).
\end{equation*}
By the differential inequality and using $k_{\max}(\omega)=+\infty$
we get \begin{equation}\label{e3.35}
\frac{1}{\omega-t}\int^{+\infty}_{k_{\max}(t)}\frac{dx}{G(x)x^{3}}\leq
1.
\end{equation}
From (\ref{e3.21}) it is easy to see that $k_{\min}(\omega)=+\infty$.
Similarly we have \begin{equation}\label{e3.36}
\frac{1}{\omega-t}\int^{+\infty}_{k_{\min}(t)}\frac{dx}{G(x)x^{3}}\geq
1.
\end{equation}
Since
\begin{equation}\label{e3.37}
\aligned
\frac{\int^{k_{\max}(t)}_{k_{\min}(t)}\frac{dx}{G(x)x^{3}}}{\int^{+\infty}_{k_{\min}(t)}\frac{dx}{G(x)x^{3}}}
&\leq\frac{\frac{k_{\max}(t)-k_{\min}(t)}{G(k_{\min}(t))k^{3}_{\min}(t)}}{\frac{1}{2G(k_{\min}(t))
k^{2}_{\min}(t)}-\frac{1}{2}\int^{+\infty}_{k_{\min}(t)}\frac{G'(x)}{G^{2}(x)x^{2}}dx}\\
&= 2\left(\frac{k_{\max}(t)}{k_{\min}(t)}-1\right)\frac{1}{1-G(k_{\min}(t))k^{2}_{\min}(t)\int^{+\infty}_{k_{\min}(t)}\frac{G'(x)}{G^{2}(x)x^{2}}dx}
.\endaligned
\end{equation}
We now claim that for the positive $z$ being so large, there holds
\begin{equation}\label{e3.38}
1-G(z) z^{2}\int^{+\infty}_{z}\frac{G'(x)}{G^{2}(x)x^{2}}dx\geq\frac{2}{2+C_{0}}.
\end{equation}
Indeed, by ($H2$) we obtain
\begin{equation*}
\int^{+\infty}_{z}\frac{G'(x)}{G^{2}(x)x^{2}}dx\leq
C_{0}\int^{+\infty}_{z}\frac{dx}{G(x)x^{3}}.
\end{equation*}
Then
\begin{equation*}
\int^{+\infty}_{z}\frac{G'(x)}{G^{2}(x)x^{2}}dx\leq
\frac{C_{0}}{C_{0}+2}\int^{+\infty}_{z}\frac{G'(x)}{G(x)x^{2}}dx+\frac{2C_{0}}{C_{0}+2}\int^{+\infty}_{z}\frac{dx}{G(x)x^{3}},
\end{equation*}
\begin{equation*}
\int^{+\infty}_{z}\frac{G'(x)}{G^{2}(x)x^{2}}dx\leq
\frac{-C_{0}}{C_{0}+2}\int^{+\infty}_{z}\left(\frac{1}{G(x)x^{2}}\right)'dx=\frac{C_{0}}{C_{0}+2}\cdot\frac{1}{G(z)z^{2}},
\end{equation*}
and this yields (\ref{e3.38}).

By (\ref{e3.32}),(\ref{e3.37}),(\ref{e3.38}) and applying $k_{\min}(\omega)=+\infty$ we have
\begin{equation*}
\lim_{t\rightarrow\omega}\frac{\int^{k_{\max}(t)}_{k_{\min}(t)}\frac{dx}{G(x)x^{3}}}{\int^{+\infty}_{k_{\min}(t)}\frac{dx}{G(x)x^{3}}}=0 .
\end{equation*}
Then we obtain
\begin{equation}\label{e3.39}
\displaystyle\lim_{t\rightarrow
\omega}\frac{\int^{+\infty}_{k_{\max}(t)}\frac{dx}{G(x)x^{3}}}{\int^{+\infty}_{k_{\min}(t)}\frac{dx}{G(x)x^{3}}}
=\lim_{t\rightarrow
\omega}\frac{\int^{+\infty}_{k_{\min}(t)}\frac{dx}{G(x)x^{3}}-\int^{k_{\max}(t)}_{k_{\min}(t)}\frac{dx}{G(x)x^{3}}}
{\int^{+\infty}_{k_{\min}(t)}\frac{dx}{G(x)x^{3}}}=1 .
\end{equation}
Combining (\ref{e3.35}), (\ref{e3.36}) with (\ref{e3.39}) we arrive at
\begin{equation*}
\lim_{t\rightarrow
\omega}\frac{1}{\omega-t}\int^{+\infty}_{k(\theta,t)}\frac{dx}{G(x)x^{3}}=1,
\quad \forall\theta\in \mathbb{S}^{1}.
\end{equation*}\hfill$\Box$

{\bf Acknowledgements.} This work is supported by the National
Natural Science Foundation of China (10671022) and Doctoral Programme
Foundation of Institute of Higher Education of China (20060027023).


\begin{thebibliography}{99}

\bibitem{BS} I.Bakas and C.Sourdis, {\it Dirichlet sigma
models and mean curvature flow},  Journal of High Energy Physics,
{\bf 6} (2007), 1088-1126.

\bibitem{GH} M.E. Gage and R.S. Hamilton, {\it The heat equation
shrinking convex plane curves},  J. Differential Geometry, {\bf 23}
(1986), 69-96.


\bibitem{B}  B.Andrews, {\it Evolving convex curves},  Calculus of
Variations and P.D.E., {\bf 7} (1998), 315-371.

\bibitem{COT} E.Calabi, P.J. Olver and A.Tannenbaum, {\it Affine
geometry, curve flows and invariant numerical approximations}, Adv.
Math, {\bf 124} (1996), 154-196.

\bibitem{ST} G.Sapiro and A.Tannenbaum, {\it On affine plane curve
evolution},  J. Functional Analysis, {\bf 119} (1994), 79-120.

\bibitem{NT}  C.H. Nien and D.H. Tsai, {\it Convex curves
moving translationally in the plane},  J. Differential equation,
{\bf 225} (2006), 605-623.

\bibitem{Ca} M.Carmo, {\it Differential Geometry of Curves
and Surfaces}, Beijing: China Machine Press, 2004.

\bibitem{Kr} N.V. Krylov, {\it Nonlinear Elliptic and Parabolic
Equations of Second Order},  D.Reidel,  1987.

\bibitem{S} R.Schneider, {\it Convex Bodies: The Brunn-Minkowski
Theory},  Cambridge University Press, 1993.

\bibitem{U} J.E. Urbas,  {\it An expansion of convex hypersurfaces},
J. Differential Geometry, {\bf 33} (1991), 91-125.

\bibitem{Z} X.P. Zhu, {\it Lectures on mean curvature flows},  AMS/IP
Studies in Advanced Mathematics, 32, American Mathematical
 Society, Providence, RI; International Press, Somerville, MA,  2002.

\bibitem{K} K.S. Chou, {\it Deforming a hypersurface by its
Gauss-Kronecker curvature}, Comm. Pure. Appl. Math, {\bf 38} (1985),
867-882.


\bibitem{B1} K.S. Chou and X.P. Zhu, {\it The curve shortening
problem}, Champman Hall/Crc (2000).

\bibitem{G1} M.E. Gage, {\it An isoperimetric inequality with
applications to curve shortening}, Duke Mathematical Journal,
{\bf50} (1984), 1225-1229.

\bibitem{G2} M.E. Gage, {\it Curve shortening makes convex curves
circular},  Invent. Math, {\bf 76} (1984), 357-364.


\bibitem{O}  R.Osserman, {\it Bonnesen-style isoperimetric
inequalities},  Amer. Math. Monthly,  {\bf 86} (1979), 1-29.









\end{thebibliography}
\end{document}